\documentclass[a4paper,10pt]{amsart}
\usepackage[english]{certus}
\usepackage[all]{xy}
\usepackage{tikz-cd}

\title{Non-amenable principal groupoids with weak containment}
\author{Vadim Alekseev}
\address{Vadim Alekseev, Technische Universit\"{a}t Dresden, Fachrichtung Mathematik, Institut f\"{u}r Geometrie, 01062, Dresden, Deutschland}
\email{vadim.alekseev@tu-dresden.de}

\author{Martin Finn--Sell}
\address{Martin Finn--Sell, Universit{\"a}t Wien, Fakult\"{a}t f\"{u}r Mathematik, Oskar-Morgenstern-Platz 1,   1090 Wien, \"{O}sterreich }
\email{martin.finn-sell@univie.ac.at}


\begin{document}

\begin{abstract}
We construct examples of principal groupoids that have weak containment but are not amenable, thus answering questions by Claire Anantharaman-Delaroche and Rufus Willett.
\end{abstract}

\maketitle

\section{Amenability and weak containment}
Through its many guises, amenability of a group has become a focal concept within both group theory and operator algebras. By a classical result of Andrzej Hulanicki (\cite{HulanickiGroupswhoseregularrepresentationweaklycontainsallunitaryrepresentations1964}), the amenability of a discrete group is equivalent to the property that all unitary representations of the group are weakly contained in the left regular representation -- we refer to this property as weak containment.

Recently, there has been interest in how far Hulanicki's classical result can be generalised, and in particular it has been shown by Rufus Willett in \cite{WillettAnonamenablegroupoidwhosemaximalandreducedCalgebrasarethesame2015} to fail for groupoids that are bundles of groups. The purpose of this note is to address Question 4.1 from \cite{Anantharaman-DelarocheSomeremarksabouttheweakcontainmentpropertyforgroupoidsandsemigroups2016} (that was also raised in Remark 3.6 of \cite{WillettAnonamenablegroupoidwhosemaximalandreducedCalgebrasarethesame2015}), namely we give an example of a \textit{principal} groupoid  that has weak containment, but is not amenable.

For a information about \'etale groupoids, their \Cs-algebras and representations, we suggest \cite[Chapter 5]{BrownCalgebrasandfinitedimensionalapproximations2008}. For more general information concerning locally compact groupoids, we refer to \cite{RenaultAgroupoidapproachtoCastalgebras1980} and \cite{Anantharaman-DelarocheAmenablegroupoids2000} and references therein.

\subsection{Preliminaries}
Throughout the text, $\mathcal{G}$ will be an \'etale Hausdorff topological groupoid, and for any subset of the unit space $U \subset \mathcal{G}^{(0)}$ we will denote by $\mathcal{G}|_{U}$ the \textit{restriction} of $\mathcal{G}$ to $U$, i.e the subgroupoid of $\mathcal{G}$ consisting of all the elements with both source and range in $U$. We remark that this groupoid is open (resp. closed) if $U$ is open (resp. closed) in $\mathcal{G}^{(0)}$.

\begin{Def}
$\mathcal{G}$ has \textit{weak containment} if the left regular representation $\lambda: C^{*}(\mathcal{G}) \rightarrow C^{*}_{r}(\mathcal{G})$ is a $*$-isomorphism. 
\end{Def}

From \cite[Theorem 6.1.4]{Anantharaman-DelarocheAmenablegroupoids2000} it is known that all measurewise amenable groupoids have weak containment. We recall the general strategy used in \cite{WillettAnonamenablegroupoidwhosemaximalandreducedCalgebrasarethesame2015} to construct a non-amenable groupoid with weak containment. 

\begin{Def}
Let $\Gamma$ be a residually finite finitely generated discrete group and let $\mc N:= \lbrace N_{i} \rbrace_{i}$ be a family of nested, finite index normal subgroups of $\Gamma$ with trivial intersection. Let $\pi_{i}$ be the quotient map $\Gamma \rightarrow \Gamma / N_{i}$. The HLS\footnote{After Nigel Higson, Vincent Lafforgue and George Skandalis who first considered this construction for a related purpose in \cite{HigsonCounterexamplestotheBaumConnesconjecture2002}.} groupoid $\mc G$ associated to $\Gamma$ and $\mc N$ is:
\begin{equation*}
\mc G:= \bigsqcup_{i \in \mathbb{N}^{+}} \lbrace i \rbrace \times X_{i}
\end{equation*}
where 
\begin{equation*}
X_{i}= \begin{cases} \Gamma/N_{i} \mbox{ if } i \in \mathbb{N}\\ \Gamma \mbox{ if } i=\infty \end{cases}
\end{equation*}
equipped with the topology generated by the following sets:
\begin{itemize}
\item the singletons $\lbrace (i,g)\rbrace$;
\item the tails: $\lbrace (i,\pi_{i}(g)) \mid i \in \mathbb{N}^{+}, i > N \rbrace$ for every fixed $g \in \Gamma$ and $N \in \mathbb{N}$.
\end{itemize}
One can check that equipped with this topology and the obvious partially defined product and inverse $\mc G$ becomes an \'etale, locally compact Hausdorff groupoid with unit space $\mathbb{N}^{+}$. Moreover, it is amenable if and only if $\Gamma$ is amenable.
\end{Def}

Considering the open invariant set $U\coloneqq \mb N \subset \mc G^{(0)}$, we obtain a commuting diagram with exact rows consisting of \Cs-algebras associated with the restriction groupoids $\mathcal{G}|_{U}$ and $\mathcal{G}|_{U^{c}}$:
\begin{equation*}
\xymatrix{
0 \ar@{->}[r] & C^{*}(\mathcal{G}|_{U})\ar@{->}[d] \ar@{->}[r]\ar@{->}[d] & C^{*}(\mathcal{G}) \ar@{->}[r] \ar@{->}[d]& C^{*}(\mathcal{G}|_{U^{c}}) \ar@{->}[r]\ar@{->}[d] & 0\\
0 \ar@{->}[r] &  C^{*}_{r}(\mathcal{G}|_{U}) \ar@{->}[r] & C^{*}_{r}(\mathcal{G}) \ar@{->}[r]  \ar@{->}[dr]& Q_{U} \ar@{->}[r] \ar@{->}[d] & 0\\
&  &  & C^{*}_{r}(\mathcal{G}|_{U^{c}}) \ar@{->}[r] & 0}
\end{equation*}
where $Q_{U}$ is the quotient by the ideal $C^{*}_{r}(\mathcal{G}|_{U})$. The groupoid $\mc G|_U$ is amenable and therefore has weak containment, and so to deduce weak containment for $\mathcal{G}$ it is enough to show that the map $C^{*}(\mathcal{G}|_{U^{c}})\to Q_{U}$ is isometric. In \cite{WillettAnonamenablegroupoidwhosemaximalandreducedCalgebrasarethesame2015}, it is then proved that $C^{*}(\mathcal{G}|_{U^{c}}) \cong C^*(\Gamma)$, $C^{*}_r(\mathcal{G}|_{U^{c}}) \cong C^*_r(\Gamma)$, and that vertical arrows come from canonical maps between these; therefore weak containment is automatic if $\Gamma$ is amenable. In the non-amenable case, weak containment is deduced from the property FD of Lubotzky--Shalom (\cite{LubotzkyFiniterepresentationsintheunitarydualandRamanujangroups2004}):

\begin{Def}
Let $\Gamma$ be a countable discrete group. $\Gamma$ has \textit{property FD} if finite dimensional representations are dense in the unitary dual of $\Gamma$; a family of finite quotients $\mc X := \lbrace \Gamma/N_{\kappa} \rbrace_{\kappa}$ is an \textit{FD family} if the set of representations of $\Gamma$ which factor through the quotient maps $\lbrace \pi_{\kappa}: \Gamma \rightarrow \Gamma / N_{\kappa} \rbrace_{\kappa}$ is dense in the unitary dual of $\Gamma$.
\end{Def}

This is then used to deduce the key result in \cite{WillettAnonamenablegroupoidwhosemaximalandreducedCalgebrasarethesame2015}: if $\mc X$ is an FD family, then the $C^{*}$-algebra $Q_{U}$ in the corresponding HLS groupoid $\mc G$ is isomorphic to the maximal group \Cs-algebra $C^{*}(\Gamma)$ through the vertical map $C^*(\Gamma)\cong C^{*}(\mathcal{G}|_{U^{c}})\to Q_{U}$ in the above diagram. However, non-amenability of $\Gamma$ implies that the groupoid $\mc G$ is non-amenable, and this finishes the construction.


\section{Constructing examples}\label{sect:groupoid}

Let $\Gamma$ be a non-amenable residually finite group with a countable nested (FD) family $\mc X$ and let $\mc G$ be the HLS groupoid from the previous section. We are going to consider a transformation groupoid constructed from $\mc G$ and the set of finite quotients $\mc X$. Let $X := \bigsqcup_{i} X_{i}$. We begin by constructing the unit space for this groupoid as a second countable compactification of $X$.

For $g\in X_i$, consider the \textit{shadow of $g$} in $X$:
\begin{equation*}
\mathop\mathrm{Sh}(g) := \bigcup_{j \geqslant i}\lbrace x \in X_{j} \mid \pi_{i,j}(x)=g \rbrace,
\end{equation*}
where $\pi_{i,j}\colon \Gamma/N_j \to \Gamma/N_i$ is the canonical quotient map.

Let $B$ be the $\mc G$-invariant \Cs-subalgebra of $\ell^\infty(X)$ generated by $\lbrace \delta_x \rbrace_{x \in X}$ and the sets of projections $\lbrace \mathbf{1}_{\mathop\mathrm{Sh}(g))} \rbrace_{g\in X_{i}}$ for all $i \in \mathbb{N}$. We will consider the spectrum of $B$, which we denote by $\widehat{X}$. As $B$ is $\mc G$-invariant, $\widehat X$ carries a natural $\mc G$-action, and we consider the transformation groupoid $G:= \widehat{X} \rtimes \mc G$. 

We remark that $G^{(0)}$ contains a obvious open invariant subset $X \subset G^{(0)}$ corresponding to the ideal generated by $\delta_x$, $x\in X$, and let $\partial X \subset G^{(0)}$ be the closed (compact) complement. The following lemma describes it as a $\Gamma$-space.

\begin{Lemma}\quad
\label{Lem:profinite}\label{lemma:direct-limit}
\begin{enumerate}
\item

 The space $\partial X$ is $\Gamma$-equivariantly homeomorphic to $\widehat{\Gamma}_{\mathcal{X}}$, the profinite completion of $\Gamma$ with respect to the family $\mathcal{X}$.
\item

 The algebra $A\coloneqq C(\partial X)$ is a direct limit of finite-dimensional $\Gamma$-\Cs-algebras $A_{i}$, such that the action on $A_i$ factors through $\Gamma/N_{i}$.
\end{enumerate}
\end{Lemma}
\begin{proof}
The inclusion $\partial X\subset \widehat X$ gives rise to a restriction homomorphism $r\colon B= C(\widehat X) \to C(\partial X) = A$ which obviously contains all elements $\delta_x$, $x\in X$, in its kernel. Thus, $A=C(\partial X)$ is generated by images of the projections $p_{i,g}\coloneqq \mathbf{1}_{\mathop\mathrm{Sh}(g)}$, $i\in\mb N,\, g\in X_i$.

Consider the finite-dimensional \Cs-algebras $A_i$ generated by the projections $p_{i,g}$, $g\in X_i$. Notice that the action of $\Gamma$ on $A_i$ obviously factors through $\Gamma/N_i$, as it is isomorphic to the natural left action of $\Gamma$ on $\mb C[\Gamma/N_i]$. Moreover, there are natural $\Gamma$-equivariant injective homomorphisms
\[
\rho_{i,j}\colon A_{i}\to A_j,
\]
\[
\rho_{i,j}(p_{i,g}) = \sum_{\pi_{i,j}(g') = g} p_{j,g'}
\]
corresponding to ($\Gamma$-equivariant) projections $\Gamma/N_j\surj \Gamma/N_i$.

Furthermore, the element $p_{i,g} - \rho_{i,i+1}(p_{i,g})$, considered as an element of $C(\widehat X)$, equals $\delta_{g}$, and therefore the kernel of the restriction map $r$ is generated by such differences. As a consequence, we get a $\Gamma$-equivariant isomorphism $A \cong \varinjlim A_i$, whence the boundary $\partial X$ is the inverse limit of the corresponding projective system of $\Gamma$-spaces. By the remark above, this projective system of spaces is naturally identified with the projective system $\{\Gamma/N_i\}_{i\in\mb N}$, equipped with the natural left $\Gamma$-action. This finishes the proof.

\end{proof}

\begin{Rem}\label{Rem:isomorphism}
As the \Cs-algebras $A_i$ are finite-dimensional, they have natural regular representations $\lambda_i\colon A_i\to \mb B(A_i)$, where $A_i$ carries the Hilbert space structure obtained from the natural trace $\tau_i\colon p_{i,g}\mapsto 1$ as well as a unitary representation $\alpha_i\colon \Gamma\to \mc U(A_i)$ given by the $\Gamma$-action. Let $\phi_{i}$ be the bijection that sends $p_{i,g}$ to $\pi_{i}(g)$. This induces an isomorphism $\phi_i$ between $A_i\rtimes \Gamma/N_i$ and the full matrix algebra $ \mb M_{|X_i|}$.
\end{Rem}

The consequence of Lemma \ref{Lem:profinite} is that we can identify the boundary piece of $G$ as 
\[
G|_{\partial X} \cong \widehat{\Gamma}_{\mc X} \rtimes \Gamma,
\]
where the latter groupoid is the transformation groupoid with the natural free action. It follows that $G$ is a principal groupoid as the action on $X$ is obviously free: $G|_X \cong \bigsqcup_{i\in \mb N} (X_i\rtimes \Gamma/N_i)$ by construction.  

Attached with this decomposition of $\widehat{X}$ into $X$ and $\partial X$ we obtain a commuting diagram with exact rows:

\begin{center}
\includegraphics{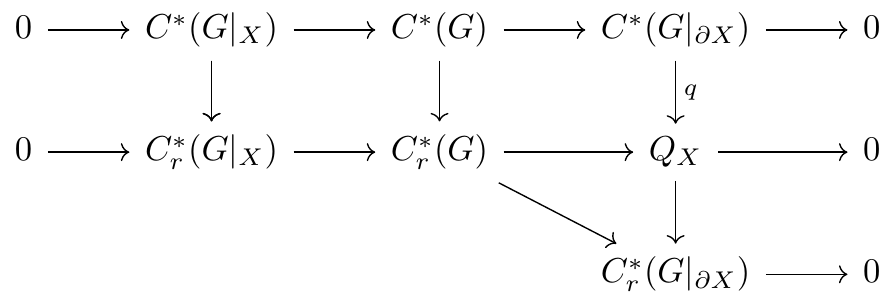}
\end{center}

\begin{Lemma}\label{lemma:weak-containment-and-isomorphisms}\quad
\begin{enumerate}
\item If the map $q\colon C^{*}({G}|_{\partial X}) \to Q_{X}$ in the above diagram is an isomorphism, then $G$ has weak containment;
\item If the map $Q_{X} \to C^{*}_{r}({G}|_{\partial X})$ in the above diagram is not an isomorphism, then $G$ is non-amenable.
\end{enumerate}
\end{Lemma}
\begin{proof}
As $G|_{X} \cong \bigsqcup_{i} X_{i}\rtimes \Gamma/N_{i}$ is a disjoint union of pair groupoids with the obvious discrete topology, it is amenable and therefore has weak containment. i) now follows from the above diagram by the five lemma. ii) follows as amenability passes to restrictions to closed invariant subsets.
\end{proof}

As a final preliminary before proving our result, we describe an ambient setting for $Q_{X}$ and $C^{*}_{r}(G)$.

\begin{Lemma}\label{lemma:embedding-into-matrices}
There is a natural isometric embedding
\[
C^*_r(G) \hookrightarrow \prod_{i\in \mb N} \mb M_{\vert X_{i}\vert},
\]
which induces an isometric embedding
\[
\iota\colon Q_{X}\hookrightarrow \frac{\prod_{j} \mb M_{\vert X_{j}\vert }}{\bigoplus_{j} \mb M_{\vert X_{j} \vert}}.
\]
\end{Lemma}
\begin{proof}
 Since $X$ is dense in $\widehat X$, we can use \cite[Corollary 2.4 a)]{KhoshkamCrossedproductsofCalgebrasbygroupoidsandinversesemigroups2004a} to see that the norm of an element $f\in C^*_r(G)$ is equal to $\sup_{x\in X} \norm{\lambda_x(f)}$, where $\lambda_x$ is the left regular representation on $s^{-1}(x)$ (which is equal to $X_i$ if $x\in X_i$). Thus we get a natural embedding
\[
C^*_r(G) \hookrightarrow \prod_{i\in \mb N} \mb M_{\vert X_{i}\vert},
\]
where $\mb M_{\vert X_{i} \vert}$ is the full matrix algebra over $X_{i}$ (viewed as bounded operators on $\ell^2(X_i)$). As $G|_X$ is a union of pair groupoids, we get $C^*_r(G|_X) \cong \bigoplus_j \mb M_{|X_j|}$, which implies that $Q_{X}$ is isometrically embedded into $\frac{\prod_{j} \mb M_{\vert X_{j}\vert }}{\bigoplus_{j} \mb M_{\vert X_{j} \vert}}$.

\end{proof}


Our goal now is connect the maximal crossed product of $A$ by $\Gamma$ with $Q_{X}$ using that $\Gamma$ has property (FD).


\begin{Prop}
The maximal crossed product $A_i\rtimes \Gamma$ embeds into $\displaystyle\prod_{j\geqslant i} A_j\rtimes \Gamma/N_j$ by the natural maps $\rho_{i,j}\rtimes \pi_j$.
\end{Prop}
\begin{proof}
The claim is equivalent to the statement that every representation of $A_i\rtimes \Gamma$ is weakly contained in a representation factoring through $A_j\rtimes \Gamma/N_j$. To this end, consider an arbitrary representation $\sigma\colon A_i \rtimes \Gamma \to \mb B(\mc H)$ and an element
\[
 x = \sum_{g\in \Gamma/N_i} p_{i,g} f_{g} \in A_{i}\rtimes_{\mathrm{alg}}\Gamma,
\]
where $f_g \in \mb C[\Gamma]$ and let $\xi,\eta\in \mc H$ be arbitrary vectors. We have
\[
 \ip{x\xi,\eta} = \sum_{g\in \Gamma/N_i} \ip{f_{g}\xi, p_{i,g}\eta}
\]

By property (FD) of $\Gamma$ for every $\eps > 0$ we get a $j\geqslant i$, representation $\sigma'\colon \Gamma \surj \Gamma/N_{j} \to \mc U(\mc H')$ and vectors $\xi'_1,\dots,\xi'_N$, $\eta'_1,\dots,\eta'_N\in \mc H'$ such that
\[
  \left|\ip{x\xi,\eta} - \sum_{g\in \Gamma/N_i} \sum_{\ell=1}^N \ip{\sigma'(f_{g}) \xi'_\ell,\eta'_\ell}\right| < \eps.
\]
Consider now the Hilbert space $\mc H'':=\mc H'\otimes A_j$ and the representation $\sigma'':=\sigma'\otimes\alpha_j \colon \Gamma\to \mc U(\mc H'')$ (which factors through $\Gamma/N_{j}$) as well as the representation $m_{i,j}\coloneqq \id_{\mc H^{'}}\otimes(\lambda_j\circ\rho_{i,j})\colon B_{i} \to \mb B(\mc H'')$. It's easy to see that these give a covariant pair and that for every $h\in \Gamma/N_j$ we have an equality of matrix coefficients
\[
 \ip{\sigma'(f_{g}) \xi'_\ell,\eta'_\ell} = \ip{\sigma''(f_{g})(\xi'_\ell\otimes p_{j,h}),\eta'_\ell\otimes \sum_{g'\in \Gamma/N_j} p_{j,g'}}
\]

Therefore any matrix coefficient of any representation of $A_{i}\rtimes \Gamma$ is approximated by a matrix coefficient of a representation factoring through $A_{j}\rtimes \Gamma/N_j$ for a suitable $j$, which ends the proof.
\end{proof}


We now can prove the following:

\begin{Prop}\label{prop:max-crossed-product-coker}
The maximal crossed product $A \rtimes \Gamma$ is isomorphic to $Q_{X}$ through the canonical quotient map $q\colon A\rtimes \Gamma \to Q_X$.
\end{Prop}
\begin{proof}
In view of Lemma \ref{lemma:embedding-into-matrices}, it is enough to prove that the composition 
\begin{equation*}
\iota\circ q\colon A\rtimes \Gamma \to \frac{\prod_{j} \mb M_{\vert X_{j}\vert}}{\bigoplus_{j} \mb M_{\vert X_{j}\vert}}
\end{equation*}
is isometric. By Lemma \ref{lemma:direct-limit} and the continuity of the maximal crossed product functor, for this it is enough to prove that the map $\iota\circ q$ is isometric on $A_i\rtimes\Gamma$.

To this end, take an arbitrary element of the algebraic crossed product $A_i\rtimes_{\mathrm{alg}}\Gamma$
\[
z = \sum_{g\in\Gamma/N_i} p_{i,g} f_g,
\]
where $f_g \in \mb C[\Gamma]$, and observe that it lifts to $C^*(G)$ as the family of elements
\[
(z_j)_j = \left(\sum_{g\in\Gamma/N_i} p_{i,g}|_{X_j} \pi_j(f_g)\right)_j \in C(X_j)\rtimes \Gamma/N_j, \quad j\geqslant i.
\]
Using the isomorphisms $\phi_j\colon A_j\rtimes \Gamma/N_j\to \mb M_{|X_j|}$ defined in Remark \ref{Rem:isomorphism}, we now see that the image of $z$ under the composition $\iota\circ q$ coincides with $(\phi_j\circ(\rho_{i,j}\rtimes \pi_j))(z)\in \prod\limits_{j\geqslant i} \mb M_{|X_j|}$, because $\rho_{i,j}(p_{i,g})(x) = p_{i,g}(x)$ for all $x\in X_j$.  


Therefore the map $\iota\circ q\colon A\rtimes \Gamma\to \frac{\prod_{j} \mb M_{\vert X_{j}\vert}}{\bigoplus_{j} \mb M_{\vert X_{j}\vert}}$ coincides with the map $A\rtimes \Gamma \to \frac{\prod_{j} \mb M_{\vert X_{j}\vert}}{\bigoplus_{j} \mb M_{\vert X_{j}\vert}}$ induced by $\phi_j\circ(\rho_{i,j}\rtimes \pi_j)$. The latter is isometric by the previous proposition, and therefore we are done.

\end{proof}

\begin{Thm}
Let $\Gamma$ be a non-amenable residually finite group with a countable nested (FD) family $\mc X$. Then the groupoid $G$ constructed above is principal and non-amenable, but has weak containment.
\end{Thm}
\begin{proof}
In view of Lemma \ref{lemma:weak-containment-and-isomorphisms} and Proposition \ref{prop:max-crossed-product-coker}, it remains to prove that the map $Q_{X} \to C^*_r(G|_{\partial X})$ is not an isomorphism. We remark that $\partial X$ has a $\Gamma$-invariant probability measure obtained by taking the weak$^*$ limit of the normalised counting measures on each $X_{i}$. By  \cite[Lemma 7.1, Remark 7.1]{WillettGeometricpropertyT2014}, we get that $Q_{X}$ contains $C^{*}(\Gamma)$ as a $*$-subalgebra, which maps onto $C^{*}_{r}(\Gamma)$ under the quotient map $Q_{X} \rightarrow C^{*}_{r}(G|_{\partial X})$. Hence the map $Q_{X} \to C^{*}_{r}(G|_{\partial X})$ is not an isomorphism. This finishes the proof.
\end{proof}

We remark that \cite[Theorems 2.2 and 2.8]{LubotzkyFiniterepresentationsintheunitarydualandRamanujangroups2004} give a wealth of examples of $\Gamma$ that satisfy the conditions above: notably free groups and surface groups (also cyclic extensions of these groups).

\section*{Acknowledgements}
This question was asked to the authors by Rufus Willett at the Erwin Schrödinger Institute programme on ``Measured group theory", February 2016. Correspondingly, the authors would like to thank both Rufus Willett for sharing the question with us and the ESI for its support. This work was also partially supported by the ERC grant “ANALYTIC” no. 259527 of Goulnara Arzhantseva.

\bibliographystyle{alpha}
\bibliography{certus-export}

\begin{thebibliography}{ADR00}

\bibitem[AD16]{Anantharaman-DelarocheSomeremarksabouttheweakcontainmentpropertyforgroupoidsandsemigroups2016}
Claire Anantharaman-Delaroche.
\newblock Some remarks about the weak containment property for groupoids and
  semigroups.
\newblock {\em arXiv:1604.01724 [math]}, April 2016.

\bibitem[ADR00]{Anantharaman-DelarocheAmenablegroupoids2000}
C.~Anantharaman-Delaroche and J.~Renault.
\newblock {\em Amenable groupoids}, volume~36 of {\em Monographies de
  L'Enseignement Math{\'e}matique [Monographs of L'Enseignement
  Math{\'e}matique]}.
\newblock {L'Enseignement Math{\'e}matique, Geneva}, 2000.
\newblock With a foreword by Georges Skandalis and Appendix B by E. Germain.

\bibitem[BO08]{BrownCalgebrasandfinitedimensionalapproximations2008}
Nathanial~P. Brown and Narutaka Ozawa.
\newblock {\em ${C}^*$-algebras and finite-dimensional approximations},
  volume~88 of {\em Graduate Studies in Mathematics}.
\newblock {American Mathematical Society}, Providence, RI, 2008.

\bibitem[HLS02]{HigsonCounterexamplestotheBaumConnesconjecture2002}
N.~Higson, V.~Lafforgue, and G.~Skandalis.
\newblock Counterexamples to the {{Baum}}\textemdash{}{{Connes}} conjecture.
\newblock {\em Geometric \& Functional Analysis GAFA}, 12(2):330--354, June
  2002.

\bibitem[Hul64]{HulanickiGroupswhoseregularrepresentationweaklycontainsallunitaryrepresentations1964}
A.~Hulanicki.
\newblock Groups whose regular representation weakly contains all unitary
  representations.
\newblock {\em Studia Math.}, 24:37--59, 1964.

\bibitem[KS04]{KhoshkamCrossedproductsofCalgebrasbygroupoidsandinversesemigroups2004a}
Mahmood Khoshkam and Georges Skandalis.
\newblock Crossed products of ${{C}^*}$-algebras by groupoids and inverse
  semigroups.
\newblock {\em Journal of Operator Theory}, 51(2):255--279, 2004.

\bibitem[LS04]{LubotzkyFiniterepresentationsintheunitarydualandRamanujangroups2004}
Alexander Lubotzky and Yehuda Shalom.
\newblock Finite representations in the unitary dual and {{Ramanujan}} groups.
\newblock In {\em Discrete geometric analysis}, volume 347 of {\em Contemp.
  Math.}, pages 173--189. {Amer. Math. Soc., Providence, RI}, 2004.

\bibitem[Ren80]{RenaultAgroupoidapproachtoCastalgebras1980}
Jean Renault.
\newblock {\em A groupoid approach to ${{C}}^{*}$-algebras}, volume 793 of {\em
  Lecture Notes in Mathematics}.
\newblock {Springer}, Berlin, 1980.

\bibitem[Wil15]{WillettAnonamenablegroupoidwhosemaximalandreducedCalgebrasarethesame2015}
Rufus Willett.
\newblock A non-amenable groupoid whose maximal and reduced ${{C}}^*$-algebras
  are the same.
\newblock {\em arXiv:1504.05615 [math]}, April 2015.

\bibitem[WY14]{WillettGeometricpropertyT2014}
Rufus Willett and Guoliang Yu.
\newblock Geometric property ({{T}}).
\newblock {\em Chinese Annals of Mathematics. Series B}, 35(5):761--800, 2014.

\end{thebibliography}

\end{document}